\documentclass[11pt]{article}
\usepackage{graphicx}
\usepackage{color}
\usepackage{amsmath}
\usepackage{amssymb}
\usepackage{amscd}
\usepackage{bbm}
\newcommand{\R}{\mathbb{R}}
\newcommand{\inr}[1]{\bigl< #1 \bigr>}

\newcommand{\E}{\mathbb{E}}

\newcommand{\PP}{{\mathbb P}}
\newcommand{\eps}{\varepsilon}

\newcommand{\cF}{{\cal F}}

\newcommand{\cC}{{\cal C}}

\newtheorem{Theorem}{Theorem}[section]

\newtheorem{Definition}[Theorem]{Definition}

\newtheorem{Corollary}[Theorem]{Corollary}
\newtheorem{Remark}[Theorem]{Remark}

\newtheorem{Assumption}{Assumption}[section]

\numberwithin{equation}{section}

\def \proof {\noindent {\bf Proof.}\ \ }

\def \endproof
{{\mbox{}\nolinebreak\hfill\rule{2mm}{2mm}\par\medbreak}}
\def\IND{\mathbbm{1}}

\begin{document}
{\baselineskip=16.5pt

\title{\bf Bounding the smallest singular value of a random matrix without concentration
}

\author{
{\bf Vladimir Koltchinskii}
\thanks{Partially supported by NSF Grants DMS-1207808, DMS-0906880 and CCF-0808863}
\\ School of Mathematics, Georgia Institute of Technology
\\ and
\\{\bf Shahar Mendelson}
\thanks{Partially supported by the Mathematical Sciences Institute -- The Australian National University and by ISF grant 900/10}
\\ Department of Mathematics, Technion -- Israel Institute of Technology
}

\maketitle

\begin{abstract}
Given $X$ a random vector in $\R^n$, set $X_1,...,X_N$ to be independent copies of $X$ and let $\Gamma=\frac{1}{\sqrt{N}}\sum_{i=1}^N \inr{X_i,\cdot}e_i$
be the matrix whose rows are $\frac{X_1}{\sqrt{N}},\dots, \frac{X_N}{\sqrt{N}}$.

We obtain sharp probabilistic lower bounds on the smallest singular value
$\lambda_{\min}(\Gamma)$ in a rather general situation, and
in particular, under the assumption that $X$ is an isotropic random vector
for which $\sup_{t\in S^{n-1}}\E|\inr{t,X}|^{2+\eta} \leq L$
for some $L,\eta>0$.
Our results imply that a Bai-Yin type lower bound holds
for $\eta>2$, and, up to a log-factor, for $\eta=2$ as well. The bounds hold without any additional assumptions on the Euclidean norm $\|X\|_{\ell_2^n}$.

Moreover, we establish a nontrivial lower bound even without any higher moment assumptions (corresponding to the case $\eta=0$), if the linear forms satisfy a weak `small ball' property.
\end{abstract}


\section{Introduction} \label{sec:intro}
The non-asymptotic theory of random matrices has attracted much attention in the recent years. One of the main questions studied in this area has to do with the behaviour of the singular values of various random matrix ensembles.

In this note, the focus is on random matrices with independent rows $({X_i}/{\sqrt{N}})_{i=1}^N,$ where $X_1,\dots, X_N$ are
independent copies of a random vector $X$ in $\R^n$; that is, if $(e_i)_{i=1}^N$ is the canonical basis in $\R^N$,
$$
\Gamma=\frac{1}{\sqrt{N}}\sum_{i=1}^N \inr{X_i,\cdot}e_i.
$$
Various attempts have been made to find conditions on $X$ that ensure that the extremal singular values of $\Gamma$ satisfy a non-asymptotic version of the classical Bai-Yin Theorem \cite{BaiYin}, formulated here.

\begin{Theorem} \label{thm:Bai-Yin}
Let $A=A_{N,n}$ be an $N \times n$ random matrix with independent
entries, distributed according to a random variable $\xi$, for
which
$$
\E \xi =0, \ \ \E \xi^2 =1 \ \ {\rm and} \ \ \E \xi^4 < \infty.
$$
If $N,n \to \infty$ and the aspect ratio $n/N$ converges to $\beta \in (0,1]$, then
$$
\frac{1}{\sqrt{N}} \lambda_{\min}(A) \to 1 -\sqrt{\beta}, \ \ \
\frac{1}{\sqrt{N}} \lambda_{\max}(A) \to 1 +\sqrt{\beta}
$$
almost surely. Also, without the fourth moment assumption, $\lambda_{\max}(A)
/\sqrt{N}$ is almost surely unbounded.
\end{Theorem}
\vskip0.3cm
In the non-asymptotic theory of random matrices, one is interested in a quantitative Bai-Yin type estimate: that for $N \geq c_1n$, with high probability, the extremal singular values of $\Gamma$ satisfy that
$$
1-c_2\sqrt{\frac{n}{N}} \leq \lambda_{\min}(\Gamma) \leq \lambda_{\max}(\Gamma) \leq 1+c_2\sqrt{\frac{n}{N}}
$$
for suitable absolute constants $c_1$ and $c_2.$

It is well understood that, if one wishes to obtain results that are non-asymptotic and that hold in the more general setup considered here (the entries of the matrix are not necessarily independent, but its rows are), one has to assume that $X$ is well behaved in some sense. Among the typical assumptions leading to a non-asymptotic Bai-Yin estimate are that $X$ is isotropic (that is, for every $t \in \R^n$, $\E\inr{X,t}^2=\|t\|_{\ell_2^n}^2$) and its distribution is a log-concave measure \cite{ALPT1,ALPT2}; that $X$ is isotropic, $\|X\|_{\ell_2^n} \leq c_1\sqrt{n}$ and $\sup_{t \in S^{n-1}} \|\inr{X,t}\|_{L_p} \leq L$ for some $p>8$ \cite{MenPao}; and that $X$ is isotropic and satisfies the assumption that, for every orthogonal projection $P$, $\|PX\|_{\ell_2^n}$ is unlikely to be large \cite{SV} (though the latter condition does not suffice to recover the $1 \pm c\sqrt{n/N}$ behaviour).

It is also known that two types of conditions are necessary for a simultaneous bound on the largest and the smallest singular values of $\Gamma$ (see, for example, the survey \cite{Ver:survey}): that with a sufficiently high probability, the Euclidean norm $\|X\|_{\ell_2^n}$ is well behaved, and that the $L_p$- and the $L_2$-norms are equivalent on linear forms for some $p>4$. These two types of assumptions -- on the Euclidean norm of $X$ and on the $L_p$ behaviour of linear forms -- are truly needed even for a weaker result: that with high probability,
$$
1-\eps \leq \lambda_{\min}(\Gamma) \leq \lambda_{\max}(\Gamma) \leq 1+\eps
$$
when $N \geq c(\eps) n$.

The main motivation for this note was the conjecture that a simultaneous estimate on the largest and smallest singular values is misleading, because their roles are very different. More precisely, that while any nontrivial estimate on the largest singular value does require concentration of some kind, usually given via an assumption on the Euclidean norm of $X$, the smallest singular value should be bounded away from zero with almost no assumptions on the random vector.

A result in this direction is due to Srivastava and Vershynin  \cite{SV}.
We state it in a slightly modified form consistent with the notations
of Theorem \ref{thm:Bai-Yin}:
\begin{Theorem} \label{thm:SV}
Let $X$ be an isotropic random vector in $\R^n$ which satisfies that for some $L,\eta>0$,
$$
\sup_{t \in S^{n-1}} \E |\inr{X,t}|^{2+\eta} \leq L.
$$
Suppose that $n/N\leq \beta$ for some $\beta\in (0,1].$
Then,
$$
\E \lambda_{\min}^2(\Gamma) \geq 1-c_1 \beta^{\frac{\eta}{2+2\eta}},
$$
for some constant $c_1>0$ depending only on $\eta, L.$
\end{Theorem}
Theorem \ref{thm:SV} implies that a slightly stronger condition than isotropicity suffices to ensure that for $N$ that is proportional to $n$, the smallest singular value is bounded away from zero, with constant probability. However, when $\eta>2$, which is the `Bai-Yin' range, or even when $\eta$ is arbitrarily large, Theorem \ref{thm:SV} does not lead to a Bai-Yin type of lower bound $1-c_1\sqrt{\beta}.$

Here, we will show that Theorem \ref{thm:SV} can be improved in the entire range of $\eta$. In particular, one has the following non-asymptotic lower bound:

\begin{Theorem} \label{thm:main-Bai-Yin}
For every $\eta>0$ and $L \geq 1$ there exist constants $c_0, c_1,\dots, c_6$  that depend only on $\eta$ and $L$ for which the following holds. Let $X$ be an isotropic vector in $\R^n$ and assume that for every $t \in S^{n-1}$ and every $u >0$,
$\PP\{|\inr{X,t}| > u\} \leq {L}/{u^{2+\eta}}$.
 
Suppose that $n/N\leq \beta$ for some $\beta\in (0,1].$ Then:
\begin{description}
\item{1.} For $\eta>2,$ with probability at least
$1-c_0 \log(e/\beta)\exp(-c_1 N\beta)$,
$$
\lambda_{\min}(\Gamma) \geq 1-c_2 \sqrt{\beta}.
$$
\item{2.} For $\eta=2,$ with probability at least
$1-\exp\Bigl(-c_3 N \beta \log (1/\beta)\Bigr)$,
$$
\lambda_{\min}(\Gamma)\geq 1 - c_4
\sqrt{\beta} \log^{3/2}\frac{1}{\beta}.
$$
\item{3.} For $0<\eta <2,$ with probability at least
$1-\exp(-c_5 N\beta \log(1/\beta))$,
$$
\lambda_{\min}(\Gamma)
 \geq 1-c_6\left(\beta
\log\left(\frac{1}{\beta}\right)\right)^{\frac{\eta}{2+\eta}}.
$$
\end{description}
\end{Theorem}

In particular, in the `Bai-Yin' range of $\eta>2$, one recovers the optimal behaviour (up to the multiplicative constants) of the smallest singular value of $\Gamma$, and with a high probability (not only
in expectation), for every $N \geq c(\beta)$.

Somewhat surprisingly, the proof of this stronger result is much easier than the proof of Theorem \ref{thm:SV} from \cite{SV}.

\vskip0.3cm

The second main result has to do with the case $\eta=0$, in which there is no additional moment information on linear forms.

It turns out that under rather minimal assumptions, one may bound the smallest singular value well away from zero.
\begin{Theorem} \label{thm:matrix-main-simple}
Let $X$ be an isotropic vector in $\R^n$ and assume that there is a constant $L$ for which $\|\inr{X,t}\|_{L_2} \leq L \|\inr{X,t}\|_{L_1}$ for every $t \in S^{n-1}$. Then, for $N \geq c_1n$, with probability at least $1-2\exp(-c_2 N)$, $\lambda_{\min}(\Gamma) \geq c_3$, for constants $c_1,c_2$ and $c_3$ that depend only on $L$.
\end{Theorem}

In fact, Theorem \ref{thm:matrix-main-simple} can be improved even further. It is reasonable to expect that there are only two situations in which the smallest singular value of $\Gamma$ is `too close' to zero. Firstly, when the distribution of $X$ is very `peaky' and almost atomic, which causes degeneracy. And secondly, when the covariance operator of $X$ is degenerate in some sense.

A more general version of Theorem \ref{thm:matrix-main-simple} shows that this belief is indeed true. The main condition used to quantify the non-degeneracy of $X$ is that the following `small-ball' function
(for one-dimensional projections of $X$)
$$
Q(u) = \inf_{t \in S^{n-1}} \PP\{|\inr{X,t}| \geq u\}
$$
is bounded away from zero {\it for some $u>0$}.

Note that $m(t,u)=\PP\{|\inr{X,t}| < u\}$ is the measure of the slab orthogonal to $t$ and of width $u$. If $X$ is absolutely continuous, then for every $t \in S^{n-1}$, $\lim_{u \to 0}m(t,u)=0$, and $m(t,u)$ is continuous in $u$ and $t$. A compactness argument shows that there is some $u>0$ for which $Q(u)>0$ (in fact, $Q(u) \geq 1/2$). Thus, the behaviour of $Q(u)$ may serve as a quantitative measure of the absolute continuity of $X$.

It turns out that having $Q(u)>0$ for some $u>0$, combined with an additional minimal condition on the covariance operator of $X,$ suffices to ensure that $\lambda_{\min}(\Gamma) \geq c$, for a constant $c$ that depends on $u$, $Q(u)$ and the covariance, with probability that grows as $1-\exp(-c_1N)$. And, the isotropicity of $X$ and the equivalence between the $L_1$ and $L_2$ norms in Theorem \ref{thm:matrix-main-simple} are only there to ensure that $Q(u)$ is bounded from below at some level, that depends only on $L$.

Theorem \ref{thm:main-Bai-Yin} and Theorem \ref{thm:matrix-main-simple} show that a sharp lower estimate on $\lambda_{\min}(\Gamma)$ is almost a `rule of nature'. It has nothing to do with a concentration phenomenon of $\|X\|_{\ell_2^n}$ and holds (almost) without any assumptions.

We end this introduction with a word about notation. Throughout this note, all absolute constants are denoted by $c,c_1,...$ or $\kappa,\kappa_1,...$. Their value may change from line to line.
We write $c(a)$ if the constant depends only on the parameter $a$; $A \lesssim B$ means that there is an absolute constant $c$ for which $A \leq cB$, and $A \lesssim_a B$ if $c$ depends only on $a$. A similar convention is used for $A \sim B$ and $A \sim_a B$, in which one has a two-sided inequality.

\section{Two main estimates} \label{sec:basic-est}

In this section, we study the problem in a more abstract setting.
Let ${\cal F}$ be a class of functions on a probability space $(\Omega,{\cal A}, P)$ and set $P_N$ to be the empirical measure $N^{-1}\sum_{i=1}^N \delta_{X_i}$, based on a sample $(X_1,\dots, X_N)$ of independent random variables with a common distribution $P$.

We will derive lower bounds on $\inf_{f\in {\cal F}}P_N f^2$
that, when ${\cal F}=\{\inr{t,\cdot}, t\in S^{n-1}\}$, imply the results on the smallest singular value of matrix $\Gamma$ stated in the introduction.

Although the estimates presented here are one-sided, they will be referred to as `isomorphic' when there is a constant $0<c<1$ for which $P_N f^2 \geq c \|f\|_{L_2(P)}^2$ for every $f \in \cF$, and as `almost isometric' when $c$ can be made arbitrarily close to $1$.

\subsection{A lower bound for a general class ${\cal F}$} \label{sec:isomorphic}

We begin with the main estimate needed to prove an
`isomorphic' bound - when $\eta=0$.

For every $u \geq 0$ let $Q_f(u)=\PP\{|f|\geq u\}$ and set
$$
Q_{\cal F}(u)=\inf_{f\in {\cal F}}Q_f(u).
$$
Let
$$
R_N({\cal F})= \E \sup_{f\in {\cal F}}\left|\frac{1}{N}\sum_{j=1}^N\eps_j f(X_j)\right|,
$$
where $(\eps_i)_{i=1}^N$ are independent, symmetric $\{-1,1\}$-valued random variables.

\begin{Theorem} \label{thm:basic}
Given a class ${\cal F}$, let $Q=Q_{\cF}$ be as above and assume that there is some $\tau>0$ for which $Q(2\tau)>0$. If
$$
R_N({\cal F})\leq \frac{\tau Q(2\tau)}{16},
$$
then
$$
{\mathbb P}\left\{
\inf_{f\in {\cal F}}\|f\|_{L_2(P_N)}^2
\geq
\frac{\tau^2 Q(2\tau)}{2}
\right\}\geq 1-2\exp\left(-\frac{Q^2(2\tau)N}{8}\right).
$$
\end{Theorem}

\proof
First, note that by Markov's inequality for the empirical measure $P_N$, $\|f\|_{L_2(P_N)}^2 \geq u^2 P_N\{|f|\geq u\}$ for every $f\in {\cal F}$ and every $u>0$. Therefore,
\begin{align*}
\|f\|_{L_2(P_N)}^2 \geq & u^2 \bigl(P\{|f|\geq 2u\} +
P_N\{|f|\geq u\}- P\{|f|\geq 2u\}\bigr)
\\
\geq & u^2 \bigl(Q(2u) +
P_N\{|f|\geq u\}- P\{|f|\geq 2u\}\bigr).
\end{align*}
Let $\phi_u:\R_+ \to [0,1]$ be defined by
\begin{equation*}
\phi_u(t) =
\begin{cases}
1 &  \ \ t \geq 2u,
\\
(t/u)-1 & \ \ u \leq t \leq 2u,
\\
0 & \ \ t<u.
\end{cases}
\end{equation*}

Observe that for every $t \in \R$, $\IND_{[u,\infty)}(t) \geq \phi_u(t)$ and $\phi_u(t) \geq \IND_{[2u,\infty)}(t)$; hence,
$$
P_N\{|f|\geq u\}\geq P_N \phi_u(|f|) \ \ {\rm and} \ \  P\{|f|\geq 2u\}\leq P\phi_u(|f|).
$$
Therefore,
$$
\|f\|_{L_2(P_N)}^2 \geq u^2 \Bigl(Q(2u) + P_N\phi_u(|f|)- P\phi_u(|f|)\Bigr).
$$
and
$$
\inf_{f\in {\cal F}}\|f\|_{L_2(P_N)}^2
\geq
u^2 \Bigl(Q(2u) -
\sup_{f\in {\cal F}}\left|P_N\phi_u(|f|)- P\phi_u(|f|)\right|\Bigr).
$$
Note that $\{\phi_u(f) : f \in F\}$ is a class of functions that is bounded by $1$. Thus, by the bounded differences inequality for
$$
H(X_1,...,X_N) = \sup_{f \in F}\left|P_N\phi_u(|f|)- P\phi_u(|f|)\right|
$$
(see, for example, \cite{BLM}) and the Gin\'{e}-Zinn symmetrization inequality \cite{GZ84}, it follows that with probability at least $1-2e^{-2t^2}$,
\begin{align*}
& \sup_{f\in {\cal F}}|P_N\phi(|f|)- P\phi(|f|)|\leq
\E \sup_{f\in {\cal F}}\left|P_N\phi(|f|)- P\phi(|f|)\right|+\frac{t}{\sqrt{N}}
\\
\leq & \frac{4}{N} \E \sup_{f \in {\cal F}} \left|\sum_{i=1}^N \eps_i \phi(|f(X_i)|) \right| + \frac{t}{\sqrt{N}}.
\end{align*}
Moreover, $\phi$ is a Lipschitz function with constant $\frac{1}{u}$; therefore, by the contraction inequality for Rademacher sums (see, e.g. \cite{LT}),
\begin{equation*}
\E \sup_{f \in {\cal F}} \left|\frac{1}{N}\sum_{i=1}^N \eps_i \phi(|f(X_i)|) \right| \leq
\frac{1}{u} \E\sup_{f\in {\cal F}}\left|\frac{1}{N}\sum_{j=1}^N \eps_j f(X_j)\right|
=
\frac{1}{u} R_N({\cal F}).
\end{equation*}
Thus, with probability at least $1-2e^{-2t^2}$,
$$
\inf_{f\in {\cal F}}\|f\|_{L_2(P_N)}^2
\geq
u^2 \left(Q(2u) - \frac{4}{u} R_N({\cal F})-\frac{t}{\sqrt{N}}
\right).
$$
Finally, if $Q(2\tau)>0$, set $u=\tau$ and put $t=\frac{Q(2\tau)\sqrt{N}}{4}$. If
$$
R_N({\cal F})\leq \frac{\tau Q(2\tau)}{16},
$$
then with probability at least $1-2\exp\left(-Q^2(2\tau)N/8\right)$,
$$
\inf_{f\in {\cal F}}\|f\|_{L_2(P_N)}^2
\geq
\frac{\tau^2 Q(2\tau)}{2},
$$
as claimed.
\endproof

If one wishes to apply Theorem \ref{thm:basic}, one has to bound $Q_{\cal F}(2\tau)$ from below. One possibility is to use the following version of the Paley-Zygmund inequality (see, e.g. \cite{GdlP}).
\begin{Theorem} \label{thm:PZ}
Let $f$ be a function on the probability space $(\Omega,P)$. For every $p,q>1$ for which $\frac{1}{p}+\frac{1}{q}=1$ and every $0<u<\|f\|_{L_1(P)}$,
$$
Q_f(u)=\PP\{|f|\geq u\} \geq \left(1-\frac{u}{\|f\|_{L_1(P)}}\right)^{q}
\left(\frac{\|f\|_{L_1(P)}}{\|f\|_{L_p(P)}}\right)^q.
$$
\end{Theorem}
\vskip0.3cm
Setting
$$
\alpha ({\cal F})=\inf_{f\in {\cal F}}\|f\|_{L_1(P)} \ \ {\rm and} \ \
\beta_p({\cal F})=\sup_{f\in {\cal F}}\frac{\|f\|_{L_p(P)}}{\|f\|_{L_1(P)}},
$$
one has the following:
\begin{Corollary} \label{cor:lower-Q}
If $\alpha({\cal F})>0$ and $\beta_p({\cal F})<\infty$ for some $p>1$, then for every $0<u<\alpha({\cal F})$,
$$
Q_{\cal F}(u)\geq
\left(1-\frac{u}{\alpha({\cal F})}\right)^{q}
\left(\frac{1}{\beta_p({\cal F})}\right)^q,
$$
where $\frac{1}{p}+\frac{1}{q} =1$.
\end{Corollary}

Theorem \ref{thm:basic} can only lead to an `isomorphic' type of lower bound, since $\tau^2 Q(2\tau)/2$ is always smaller than $\inf_{f \in F} \|f\|_{L_2(P)}^2$. However, we will show in the next section that in some cases it is possible to obtain an `almost isometric' bound as well.

\subsection{A sharp lower bound} \label{sec:almost-isometric-vc}
Here, we will prove a high probability `almost isometric' lower bound, of the form $\inf_{f \in \cF} P_N f^2 \geq 1-\delta$, for certain subsets of the unit sphere in $L_2(P)$.

\begin{Definition} \label{def:VC}
Let $\cC$ be a class ${\cal C}$ of subsets of $\Omega$. A set $\{x_1,...,x_n\}$ is shattered by $\cC$ if for every $I \subset \{1,...,n\}$ there is some $C \in \cC$ for which $x_i \in C$ if $i \in I$, and $x_i \not \in C$ otherwise.

The maximal cardinality of a subset of $\Omega$ that is shattered by $\cC$ is called the VC dimension of $\cC$, and is denoted by $VC(\cC)$.
\end{Definition}
We refer the reader to \cite{VW} for basic facts on VC classes and on the VC-dimension.

\begin{Assumption} \label{ass:almost-isometric-classes}
Let $\cF$ be a subset of the $L_2(P)$ unit sphere. Assume that
\begin{description}
\item{1.} there is a constant $L>1$ for which, for every $u >0$,
\begin{equation} \label{eq:tail-estimate}
\sup_{f \in \cF} \PP\{|f| > u\} \leq \frac{L}{u^{2+\eta}}.
\end{equation}
\item{2.} The collection of sets ${\cC}=\{\{|f|> u\}: f\in {\cF}, u>0\}$ is a VC class of sets of VC-dimension at most $d$.
\end{description}
\end{Assumption}

\vskip0.3cm

Clearly, the first part of Assumption \ref{ass:almost-isometric-classes} follows if the $L_{2+\eta}$ norm is equivalent to the $L_2$ one on $\cF$. And, as will be explained later, the class of norm-one linear functionals on $\R^n$ satisfies the second part for $d \sim n$.

\begin{Theorem} \label{thm:emp-iso}
For every $\eta>0$ and $L \geq 1,$ there exist constants $\kappa_i$, $i=0,1,...,6$ that depend only on $L$ and $\eta$ for which the following holds. Assume that $\cF$ satisfies Assumption \ref{ass:almost-isometric-classes} with constants $\eta$, $L$ and $d$, and that
$d/N\leq \beta$ for some $\beta \in (0,1]$. Then,
the following statements hold:
\begin{description}
\item{1.} If $\eta>2,$ then with probability at least
$1-\kappa_0\log (e/\beta)\exp(-\kappa_1 N\beta),$
$$
\inf_{f \in \cF} P_N f^2 \geq 1 - \kappa_2 \sqrt{\beta}.
$$
\item{2.} If $\eta=2,$ then with probability at least
$1-\exp\Bigl(-\kappa_3 N \beta \log(1/\beta)\Bigr),$
$$
\inf_{f \in \cF} P_N f^2 \geq 1 - \kappa_4
\sqrt{\beta}\log^{3/2}\frac{1}{\beta}.
$$
\item{3.} If $0<\eta <2,$ then with probability at least
$1-\exp\Bigl(-\kappa_5 N\beta \log(1/\beta)\Bigr)$,
$$
\inf_{f \in \cF} P_N f^2 \geq 1 - \kappa_6 \left(\beta \log \left(\frac{1}{\beta}\right)\right)^{\frac{\eta}{2+\eta}}.
$$
\end{description}
\end{Theorem}

We will present a detailed proof of the first part of Theorem \ref{thm:emp-iso}. Since the other two parts are almost identical to the first one, we will only outline the minor differences in their proofs.

The following well-known fact regarding VC classes has a key role in the proof of Theorem \ref{thm:emp-iso}.
\begin{Theorem} \label{thm:VC}
There exists an absolute constant $\kappa$ for which the following holds.
Let $\cC$ be a class of sets, put $d={\rm VC}({\cal C})$ and set $\sigma^2 = \sup_{C \in \cC} P(C)$. For every $t>0$, with probability at least $1-2\exp(-t)$,
$$
\sup_{C \in \cC}|P_N(C)-P(C)| \leq \kappa \left(\sigma \sqrt{\frac{d}{N} \log\left(\frac{e}{\sigma}\right)}+\frac{d}{N}\log\left(\frac{e}{\sigma}\right)+\sigma \sqrt{\frac{t}{N}}+\frac{t}{N} \right).
$$
\end{Theorem}
Theorem \ref{thm:VC} follows from Talagrand's concentration inequality for uniformly bounded
empirical processes (see, e.g., \cite{Led-book}) and a standard estimate on $\E\sup_{C \in \cC}|P_N(C)-P(C)|$ \cite{Tal94}.

\vskip0.5cm
\noindent{\bf Proof of Theorem \ref{thm:emp-iso}.} For $j \geq 1$, let $\cC_j=\{ \{|f| > u\} : f \in \cF, \ 2^j \leq u \leq 2^{j+1}\}$ and denote  $\cC_0=\{ \{|f| > u\} \ : \ f \in \cF, \ 0 < u \leq 1\}$.

Recall that by Assumption \ref{ass:almost-isometric-classes}, $VC(\cC) \leq d$. Since $\cC_j \subset \cC$, it follows that ${\rm VC}(\cC_j) \leq d$ for $j=0,1,..$. Moreover, setting $\sigma_j = \sup_{C \in \cC_j} (P^{1/2}(C))$, it is evident from the tail estimate in Assumption \ref{ass:almost-isometric-classes}  that for $j \geq 1$,
$$
\sigma_j^2 =\sup_{f \in \cF} \PP\{|f| \geq 2^j\} \leq L2^{-j(2+\eta)},
$$
and trivially, that $\sigma_0 \leq 1$.
Therefore, by Theorem \ref{thm:VC}, for every $j \geq 1$, with probability at least $1-2\exp(-t)$,
\begin{align} \label{eq:est-on-sets}
& \sup_{C \in \cC_j} |P_N(C)-P(C)| \nonumber
\\
\leq & \kappa_0 \left(\sigma_j \sqrt{\frac{d}{N}\log\left(\frac{e}{\sigma_j}\right)}+\frac{d}{N}\log\left(\frac{e}{\sigma_j}\right)+\sigma_j \sqrt{\frac{t}{N}}+\frac{t}{N} \right) \nonumber
\\
\leq & \kappa_1 L^{1/2}(1+\eta) \left( 2^{-j(1+\eta/2)}\left(\sqrt{j\frac{d}{N}}+\sqrt{\frac{t}{N}}\right)+ j\frac{d}{N}+\frac{t}{N} \right),
\end{align}
for suitable absolute constants $\kappa_0$ and $\kappa_1$. A similar estimate holds for $j=0$.

Fix $\delta$ to be specified later and set $A=\max\{(L/\eta \delta)^{1/\eta},1\}$. We will assume that $A>1$, as the case $A=1$ is considerably simpler and is omitted.

Note that for every $f \in \cF$,
$$
2\int_A^\infty uP\{|f|>u\}du \leq 2 \int_A^\infty L u^{-(1+\eta)}du \leq \delta.
$$
Recall that $\|f\|_{L_2(P)}=1$ for every $f \in \cF$, and thus
\begin{align*}
P_N f^2  \geq & 2\int_0^A uP_N\{|f|>u\}du
\\
= & 2\int_0^A uP\{|f|>u\}du +2\int_0^A u\left(P_N\{|f|>u\}-P\{|f|>u\}\right)du
\\
\geq & 1-\delta - 2\int_0^A u\big|P_N\{|f|>u\}-P\{|f|>u\}\big| du.
\end{align*}
Observe that for every $f \in \cF$ and $u \in [2^j,2^{j+1}]$,
$$
\left|P_N\{|f|>u\}-P\{|f|>u\}\right| \leq \sup_{C \in \cC_j} |P_N(C)-P(C)|.
$$
Hence, if $j_0$ is the smallest integer for which $2^{j_0} \geq A$,
\begin{align*}
& \int_0^A u\left|P_N\{|f|>u\}-P\{|f|>u\}\right| du
\\
\leq & \int_0^1 u\left|P_N\{|f|>u\}-P\{|f|>u\}\right| du
\\
 + & \sum_{j=1}^{j_0} \int_{2^j}^{2^{j+1}} u\left|P_N\{|f|>u\}-P\{|f|>u\}\right| du
\\
\lesssim &
\sum_{j=0}^{j_0} 2^{2j} \sup_{C \in \cC_j} |P_N(C)-P(C)| =(*).
\end{align*}
Applying \eqref{eq:est-on-sets} and summing the probabilities, it is evident that with probability at least $1-2(j_0+1)\exp(-t)$, for every $0 \leq j \leq j_0$
\begin{align}
\label{terms}
& 2^{2j} \sup_{C \in \cC_j} |P_N(C)-P(C)|
\\
\nonumber
\lesssim &
L^{1/2}(1+\eta)  2^{j(1-\eta/2)}\left(\sqrt{j\frac{d}{N}}+\sqrt{\frac{t}{N}}\right)+2^{2j} \left(j\frac{d}{N}+\frac{t}{N}\right).
\end{align}

\noindent 1. If $\eta>2$, then (\ref{terms}) implies that
\begin{equation}
\label{bound*}
(*) \lesssim_{\eta,L} \biggl(\sqrt{\frac{d}{N}}+\sqrt{\frac{t}{N}}
+  \frac{d}{N} A^2 \log A + A^2 \frac{t}{N}\biggr)  = (**).
\end{equation}
Under the assumption that $d/N\leq \beta$ for $\beta \in (0,1],$
set $t=\kappa_2 N\beta$ for a constant $\kappa_2$ to be selected later, and that depends only on $L$ and on $\eta$, and let $\delta=\sqrt{\beta}$; hence,
$$
A\lesssim_{\eta,L} \beta^{-1/(2\eta)}.
$$
By a straightforward computation, combined with a proper choice
of the constant $\kappa_2$ in the definition of $t$,
$$
(**) \lesssim_{\eta,L}\sqrt{\beta}.
$$
To estimate the probability of that event, note that
$$
2^{j_0} \lesssim A\lesssim_{\eta,L} \beta^{-1/(2\eta)}
$$
implying that
$$
j_0+1 \lesssim_{\eta, L} \log(e/\beta).
$$
Thus, $2(j_0+1) \exp(-t) \leq \kappa_3 \log (e/\beta)\exp(-\kappa_2 N \beta)$, and
with probability at least
$1-\kappa_3 \log (e/\beta)\exp\left(-\kappa_2 N\beta\right)$,
\begin{align*}
\inf_{f \in \cF} P_N f^2 \geq 1-\kappa_4\sqrt{\beta},
\end{align*}
completing the proof of the first part.

\noindent 2. In the case $\eta=2,$ (\ref{bound*}) becomes
$$
(*) \lesssim_{\eta,L} \left((\log A)^{3/2}\sqrt{\frac{d}{N}}+\log A\sqrt{\frac{t}{N}}
+  \frac{d}{N} A^2 \log A + A^2 \frac{t}{N}\right)  = (**).
$$
One may take $\delta= \sqrt{\beta}\log^{3/2}(1/\beta)$ and $t= \kappa_5 N\beta \log (1/\beta)$ and repeat the argument used in the case $\eta>2$.

\noindent 3. Finally, if $\eta<2$, then (\ref{bound*}) turns out to be
\begin{equation*}
(*) \lesssim_{\eta,L}  A^{1-\eta/2}\sqrt{\log A}\sqrt{\frac{d}{N}}+A^{1-\eta/2}\sqrt{\frac{t}{N}}
+  \frac{d}{N} A^2 \log A + A^2 \frac{t}{N} = (**),
\end{equation*}
and one can again repeat the same argument used in the case
$\eta>2$, this time for the choice
$$
t=\kappa_6 N\beta \log(1/\beta) \ \ {\rm  and} \ \  \delta=\left(\beta \log(1/\beta)\right)^{\eta/(2+\eta)}
$$
for a constant $\kappa_6$ that depends only on $L$ and $\eta$.
\endproof

\section{The smallest singular value of a random matrix} \label{sec:smallest}
Let $X$ be a random vector on $\R^n$, set $X_1,...,X_N$ to be independent copies of $X$ and put
$\Gamma=\frac{1}{\sqrt{N}} \sum_{i=1}^N \inr{X_i,\cdot}e_i$ -- the corresponding random matrix. Clearly, the smallest singular value of $\Gamma$ satisfies
$$
\lambda_{\min}(\Gamma) = \inf_{t \in S^{n-1}} \|\Gamma t\|_{\ell_2^N} = \inf_{t \in S^{n-1}} \|\inr{t,\cdot}\|_{L_2(P_N)}.
$$
Thus, it seems natural to apply Theorem \ref{thm:basic} and Theorem \ref{thm:emp-iso} for the choices of ${\cal F}=\{\inr{t,\cdot} : t \in S^{n-1}\}$ and $P$ -- the probability distribution of random vector $X$.
\vskip0.3cm
\noindent{\bf Proof of Theorem \ref{thm:main-Bai-Yin}.} It suffices to show that Assumption \ref{ass:almost-isometric-classes} holds for the class $\cF=\{\inr{t,\cdot} : t \in S^{n-1}\}$ and apply Theorem \ref{thm:emp-iso}. By Markov's inequality, individual tail estimates of the form $\PP\{|\inr{X,t}|>u\} \leq L/u^{2+\eta}$ are true when the $L_2$ and $L_{2+\eta}$ norms are equivalent on $\R^n$, as it is assumed here. To complete the proof, one only has to verify that the class of sets
$$
\cC=\left\{ \left\{|\inr{t,\cdot}| >u\right\}, \ t \in S^{n-1}, \ u>0\right\}
$$
satisfies $VC(\cC) \lesssim n$.

It is well known that the class of halfspaces in $\R^n$ is a VC class of dimension at most $n+1$ (see, e.g., \cite{VW}). It is also standard to verify that if ${\cal A}$ and ${\cal B}$ are VC class with $\max\{ VC({\cal A}),VC({\cal B})\} \leq d,$ then for ${\cal D}=\{A \cup B : A \in {\cal A}, B \in {\cal B}\}$, $VC({\cal D}) \lesssim d$.

Here, $\cC$ is contained in a class of sets of that form, with ${\cal A}$ and ${\cal B}$ being the sets of halfspaces in $\R^n$. Thus, $VC(\cC) \leq \kappa n$ for some absolute constant $\kappa$ and Theorem \ref{thm:main-Bai-Yin} follows from Theorem \ref{thm:emp-iso}.
\endproof

Next, let us consider the case $\eta=0$. The next result shows that the smallest singular value of $\Gamma$ is bounded from below in a more general setup than the one formulated in Theorem \ref{thm:matrix-main-simple}.

\begin{Assumption} \label{ass:random-vector}
Let $X$ be a random vector on $\R^n$ and assume that
\begin{description}
\item{1.} There exist constants $0<a<A$ for which $a \leq \|\inr{X,t}\|_{L_2} \leq A$ for every $t \in S^{n-1}$.
\item{2.} There exists a constant $B$ that satisfies that $\|\inr{X,t}\|_{L_2} \leq B \|\inr{X,t}\|_{L_1}$ for every $t \in S^{n-1}$.
\end{description}
\end{Assumption}
The first part of Assumption \ref{ass:random-vector} simply states that the covariance operator of $X$ is non-degenerate. The second is the additional component that allows one to bound the `small-ball' function $Q$ from below.

\begin{Theorem} \label{thm:main-isomorph-matrix}
There exist absolute constants $c_0$, $c_1$ and $c_2$ for which the following holds. If $X$ satisfies Assumption \ref{ass:random-vector} and $N \geq c_0 B^4 (A/a)^2 n$, then with probability at least $1-\exp(-c_1B^4N)$, $\lambda_{\min}(\Gamma) \geq c_2a/B^2$.
\end{Theorem}

\proof Let ${\cal F}=\left\{\inr{t,\cdot} : t \in S^{n-1}\right\}$ and let $P$ be the distribution of $X.$ Using the notation of Corollary \ref{cor:lower-Q},
\begin{equation*}
\alpha({\cal F})= \inf_{t \in S^{n-1}} \|\inr{t,X}\|_{L_1} = \inf_{t \in S^{n-1}} \|\inr{t,X}\|_{L_2} \cdot \frac{\|\inr{t,X}\|_{L_1}}{\|\inr{t,X}\|_{L_2}}
\geq  \frac{a}{B},
\end{equation*}
and
\begin{equation*}
\beta_2({\cal F}) = \sup_{t \in S^{n-1}} \frac{\|\inr{t,X}\|_{L_2}}{\|\inr{t,X}\|_{L_1}} \leq B.
\end{equation*}
Therefore, for every $0<u<a/B$,
$$
Q_{\cal F}(u)\geq \left(1-\frac{B u}{a}\right)^2 \left(\frac{1}{B}\right)^2,
$$
and setting $\tau= a/4B$, it follows that $Q(2\tau)\geq {1}/{4B^2}$.

To bound $R_N({\cal F})$, observe that

$$
R_N({\cal F})=\E \sup_{t\in S^{n-1}}\left|\frac{1}{N}\sum_{j=1}^N \eps_j \inr{t,X_j}\right|=
\E \left\|\frac{1}{N}\sum_{j=1}^N \eps_j X_j\right\|_{\ell_2^n},
$$
and since $\E\|X\|_{\ell_2^n}^2= \sum_{i=1}^n \E \inr{X,e_i}^2 \leq A^2n$, it is evident that
\begin{equation*}
\E\left\|\frac{1}{N}\sum_{j=1}^N \eps_j X_j\right\|_{\ell_2^n}\leq
\left(\E \left\|\frac{1}{N}\sum_{j=1}^N \eps_j X_j\right\|_{\ell_2^n}^2\right)^{1/2}
=
\sqrt{\frac{\E\|X\|_{\ell_2^n}^2}{N}} \leq A\sqrt{\frac{n}{N}}.
\end{equation*}
These estimates, combined with Theorem \ref{thm:basic}, show that if $R_N({\cal F}) \leq \tau Q(2\tau)/16$, there are absolute constants $c_1$ and $c_2$, for which, with probability at least $1-\exp(-c_1N/B^4)$,
$$
\lambda_{\min}(\Gamma) =\inf_{t \in S^{n-1}} \|\Gamma t\|_{\ell_2^N} \geq c_2 \frac{a}{B^2}.
$$
This is the case when
$$
A \sqrt{\frac{n}{N}} \leq \frac{a}{(16B)^2}.
$$
\endproof

\begin{Remark} \label{rem:singular-value-general}
Observe that even weaker assumption than Assumption \ref{ass:random-vector} suffices to bound the smallest singular value of $\Gamma$. The proof shows that if $(\E \|X\|_{\ell_2^n}^2)^{1/2} \leq A \sqrt{n}$, $Q(2\tau)>0$ for some $\tau$ and $N \geq c_1 An/\tau^2 Q^2(2\tau)$, then $\lambda_{\min}(\Gamma) \geq c_2 \tau Q^{1/2}(2\tau)$ with probability $1-2\exp(-c_3NQ(2\tau)^2)$, for absolute constants $c_1,c_2$ and $c_3$.
\end{Remark}

\footnotesize {

\end{document}
\begin{thebibliography}{10} \frenchspacing
\bibitem{ALPT1} R. Adamczak, A. Litvak, A. Pajor, N.
Tomczak-Jaegermann, Quantitative estimates of the convergence of
the empirical covariance matrix in log-concave ensembles, Journal of the American
Mathematical Society 23 535-561, 2010.

\bibitem{ALPT2} R. Adamczak, A. Litvak, A. Pajor, N.
Tomczak-Jaegermann, Sharp bounds on the rate of convergence of
empirical covariance matrix, C.R. Math. Acad. Sci. Paris, 349,
195--200, 2011.

\bibitem{BaiYin} Z.D. Bai, Y.Q. Yin, Limit of the smallest eigenvalue of a large dimensional sample covariance matrix, Ann. Probab. 21, 1275--1294, 1993.


\bibitem{BLM} S. Boucheron, G. Lugosi, and P. Massart, {\it Concentration Inequalities: A Nonasymptotic Theory of Independence}, Oxford University Press, 2013.

\bibitem{GdlP}  V. de la Pe\~{n}a, E. Gin\'{e} {\it Decoupling: From Dependence to Independence}, Springer-Verlag, 1999.

 \bibitem{GZ84} E. Gin\'{e}, J. Zinn, Some limit theorems for empirical
processes, Annals of Probability 12(4), 929-989, 1984.


\bibitem{Led-book} M. Ledoux, {\it The Concentration of Measure Phenomenon}, Mathematical Surveys and Monographs 89, AMS, 2001.

\bibitem{LT} M. Ledoux, M. Talagrand,  {\it Probability in
  Banach spaces. Isoperimetry and processes}, Ergebnisse der Mathematik
  und ihrer Grenzgebiete (3), vol. 23.  Springer-Verlag, Berlin, 1991.

\bibitem{MenPao} S. Mendelson, G. Paouris, On the singular values of random matrices, Journal of the European Mathematics Society, to appear.

\bibitem{SV} N. Srivastava, R. Vershynin, Covariance estimation for distributions with $2+\epsilon $ moments, Annals of Probability 41 (2013), 3081--3111.

\bibitem{Tal94} M. Talagrand, Sharper Bounds for Gaussian and Empirical Processes, Annals of Probability, 22(1) 28-76, 1994.

\bibitem{VW}{A.W. Van der Vaart, J.A. Wellner, }
{\it Weak convergence and empirical processes}, Springer Verlag,
1996.

\bibitem{Ver:survey} R. Vershynin, Introduction to the non-asymptotic analysis of random matrices. In: {\it Compressed Sensing: Theory and Applications}, Yonina Eldar and Gitta Kutyniok (eds), 210-268, 2012, Cambridge University Press.


\end{thebibliography}
